\documentclass[12pt]{article}

\usepackage{amsthm}
\usepackage{amsmath}
\usepackage{amsfonts}
\usepackage{amssymb}
\usepackage{fullpage}
\usepackage{stmaryrd}
\usepackage{graphicx}
\usepackage{psfrag}

\usepackage{pslatex}
\newtheorem{lemma}{Lemma}
\newtheorem{cor}[lemma]{Corollary}
\newtheorem{theorem}[lemma]{Theorem}

\theoremstyle{definition}
\newtheorem{defn}[lemma]{Definition}

\newcommand{\Ref}[1]{(\ref{#1})}

\newcommand{\Ddiff}[3]{ \frac{\partial^{#3} \! #1}{\partial #2^{#3}}} 

\newcommand{\nn}{\nonumber \\}

\newcommand{\spacebreak}{\begin{displaymath}
  \triangleleft \; \lhd \;
  \diamond \;
  \rhd \; \triangleright
\end{displaymath}}

\newcommand{\nds}{\preceq_s}

\newcommand{\seq}[1]{\llbracket \, #1 \, \rrbracket}

\newcommand{\hhp}[1]{{| #1|_{\Leftrightarrow}}}
\newcommand{\vhp}[1]{{| #1|_{\Updownarrow}}}

\newcommand{\hp}[1]{ \, {\odot}_{#1} \, }

\begin{document}

\author{Andrew Rechnitzer  \\ 
  {\small Department of Mathematics and Statistics,} \\
  {\small The University of Melbourne, Parkville Victoria 3010, Australia.}\\
  {\small email: \texttt{A.Rechnitzer@ms.unimelb.edu.au}} \\
}
\title{Haruspicy 3: \\ The directed bond-animal generating function \\ is not D-finite.}
\maketitle

\begin{abstract}
  While directed site-animals have been solved on several lattices,
  directed bond-animals remain unsolved on any non-trivial lattice. In
  this paper we demonstrate that the anisotropic generating function
  of directed bond-animals on the square lattice is fundamentally
  different from that of directed site-animals in that it is not
  differentiably finite. We also extend this result to directed
  bond-animals on hypercubic lattices.
  
  This indicates that directed bond-animals are unlikely to be solved
  by similar methods to those used in the solution of directed
  site-animals. It also implies that a solution cannot by conjectured
  using computer packages such as GFUN \cite{Algolib} or differential
  approximants \cite{Guttmann1989}.
\end{abstract}

\section{Introduction}
The enumeration of lattice animals is a long-standing problem in
enumerative combinatorics and finds applications in statistical
physics and theoretical chemistry. Though the subject has received
considerable attention over many years, the problem remains unsolved.

\begin{figure}[ht]
  \centering
  \includegraphics[scale=0.5]{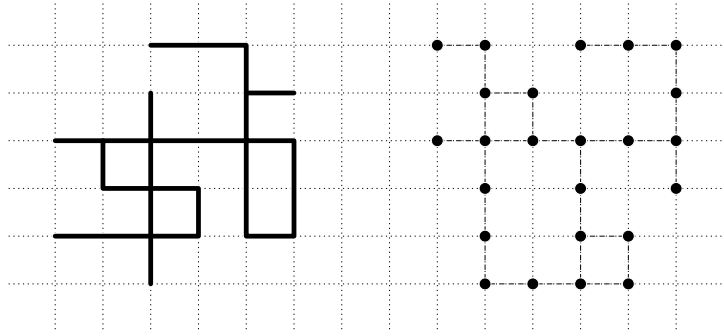}
  \caption{A bond-animal and a site-animal.}
  \label{fig animal eg}
\end{figure}

\begin{defn}
  A \emph{bond-animal} is a connected union of bonds (edges) on a
  lattice\footnote{Except for the animals in Corollary~\ref{cor dfinite hypercube}, 
    all objects considered in this paper are on the square lattice.}.
  Similarly a \emph{site-animal} is a connected union of sites
  (vertices). Two animals are considered to be the same if they are
  translates of each other.
\end{defn}

In spite of the difficulty of enumerating general lattice animals,
many subclasses have been solved. In almost all cases it has only been
possible to count animals with quite severe topological restrictions
--- such as directedness or convexity. In this paper we focus on
directed animals.

\begin{figure}[ht]
  \centering
  \includegraphics[scale=0.5]{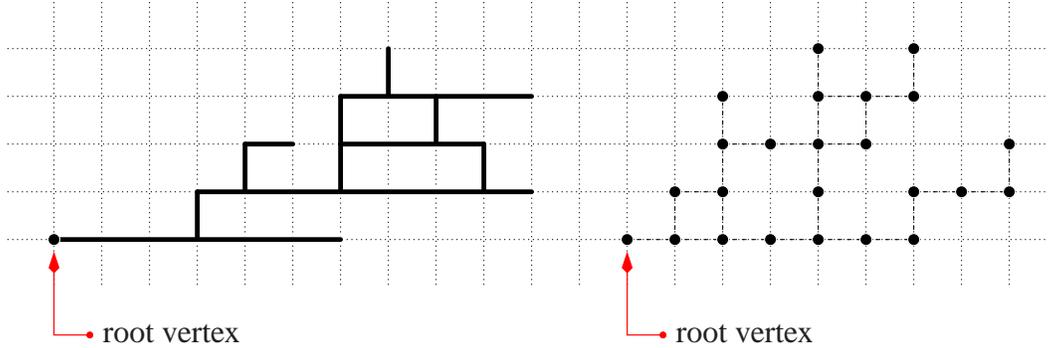}
  \caption{A directed bond-animal and a directed site-animal.}
  \label{fig dba eg}
\end{figure}

\begin{defn}
  A bond-animal is \emph{directed} if it contains a
  special vertex called the \emph{root vertex} such that all bonds in
  the animal may be reached from the root vertex by paths taking only
  north and east steps. Similarly a site-animal is directed if
  it contains a root vertex and all other sites can be reached from it
  by taking only north and east steps. See Figure~\ref{fig dba eg}
\end{defn}

Directed site-animals were first solved around 20 years ago by Dhar
\cite{dhar-premiere, dhar-three} by mapping the problem to a hard-core
lattice gas, and then subsequently by a number of authors using more
geometric and bijective methods (such as
\cite{betrema-penaud-guingois, betrema-penaud, gouyou-viennot,
  penaud}). The resulting generating function is a simple algebraic
function:
\begin{equation}
  \label{eqn dsa gf}
  S(q) = \sum_{A \in \, \substack{\mathrm{directed} \\ \mathrm{site-animals}}} q^{|A|} = 
  \frac{1}{2} \left(  \sqrt{ \frac{1+q}{1-3q} } -1 \right),
\end{equation}
where $|A|$ denotes the number of sites in an animal, $A$.  A similar
solution exists for directed site-animals on the triangular lattice
and a directed cubic lattice (in which both nearest-neighbour and
next-nearest-neighbour steps are allowed).

The generating function of directed bond-animals is defined in a
similar way:
\begin{equation}
  B(z) = \sum_{A \in\, \substack{\mathrm{directed} \\ \mathrm{bond-animals}} } z^{|A|},
\end{equation}
where $|A|$ denotes the number of bonds in the animal $A$. Despite the
similarity of the underlying objects, the directed bond animal
generating function remains unsolved.

In this paper we show that a possible reason that directed
bond-animals remain unsolved is that their generating function, in
particular their \emph{anisotropic} generating function, is not within
the class of \emph{differentiably finite} functions. Consequently it
is is fundamentally different from that of directed site-animals and
most other solved bond lattice models. A similar result for
self-avoiding polygons was recently given in \cite{ADR_sap}.

In the next section we define differentiably finite functions and the
anisotropic generating functions of directed bond and directed site
animals. In Section~\ref{sec dfinite proof} we prove that the
anisotropic generating function of directed bond-animals is not
differentiably finite. An immediate corollary of this is that the
generating function of directed bond-animals on the $d$-dimensional
hypercubic lattice (with $d \geq 2$) is not D-finite.

\section{Anisotropic and differentiably finite generating functions}
Perhaps the most common functions in combinatorics and mathematical
physics are those that satisfy simple linear differential equations
with polynomial coefficients --- these functions are called
\emph{differentiably finite} or \emph{D-finite}. More precisely:

\begin{defn}
  Let $f(t)$ be a formal power series in $t$ with coefficients in
  $\mathbb{C}$. This series is \emph{differentiably finite} or
  \emph{D-finite} if there exist a nontrivial differential equation of
  the form
  \begin{equation}
    \label{eqn dfinite defn}
    P_k(t) f^{(k)}(t) + \dots + P_1(t) f'(t) + P_0(t) f(t) = 0,
  \end{equation}
  where the $P_i(t)$ are polynomials in $t$ with complex coefficients.
  It can also be shown that any algebraic power series is also a
  D-finite power series \cite{Lipshitz1989}.
\end{defn}

Ideally we would like to show that the generating function, $B(z)$, is
fundamentally different in nature from that of directed site-animals,
$S(q)$, which is an algebraic, and hence D-finite, power series.
Perhaps the easiest way to demonstrate that a series is not D-finite
is to examine its singularities; the classical theory of linear
differential equations implies that D-finite series of a single
variable cannot have an infinite number of singularities. By this
reasoning the function $f(t) = \tan(t)$ is not a D-finite power
series in $t$.

Unfortunately, almost nothing is known rigorously about $B(z)$ --- we
do not even know the exact location of its dominant singularity, and
(the author) certainly cannot show that it has an infinite number of
singularities. Fortunately, by considering the \emph{anisotropic
  generating function} we are able to make considerably more progress.

\begin{figure}[ht]
  \centering
  \includegraphics[scale=0.5]{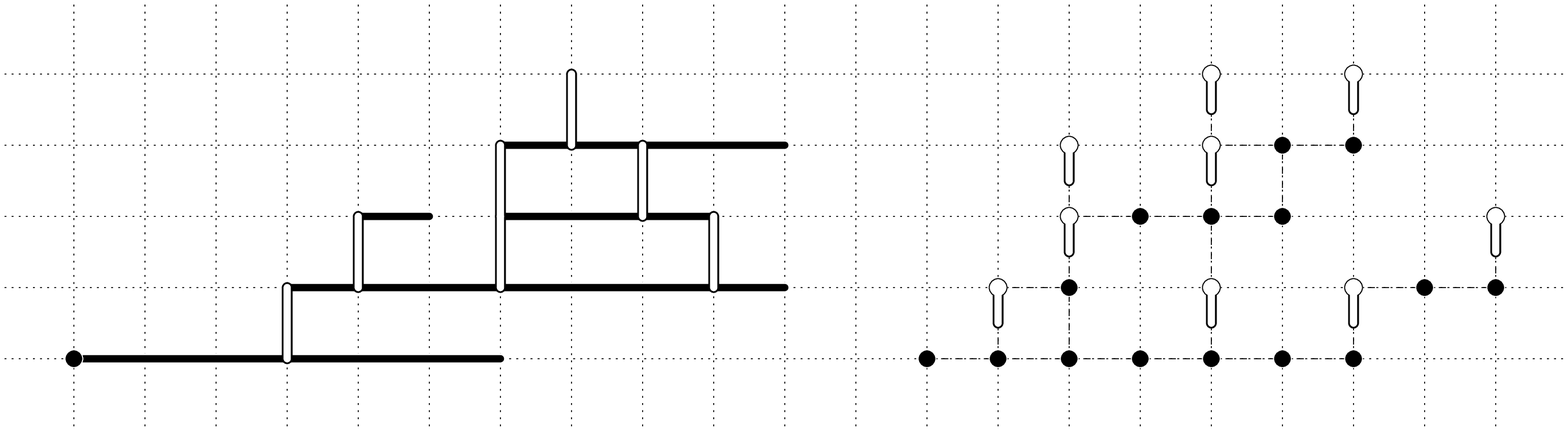}
  \caption{(left) A directed bond-animal with vertical bonds
    highlighted. (right) A directed site-animal with sites supported
    only from the south highlighted.}
  \label{fig dba-aniso eg}
\end{figure}

We form the anisotropic generating function of directed bond-animals
by counting animals, $A$, according to the number of horizontal
(resp. vertical) bonds it contains, denoted $\hhp{A}$ (resp. $\vhp{A}$):
\begin{equation}
  B(x,y) = \sum_{A \in \, \substack{\mathrm{directed} \\ \mathrm{bond-animals}}}
  x^{\hhp{A}} y^{\vhp{A}}
  = \sum_{m,n} b_{m,n} x^m y^n,
\end{equation}
where $b_{m,n}$ is the number of directed bond-animals with $m$
horizontal and $n$ vertical bonds.

Anisotropisation of the directed site-animal generating function is
more problematic and is not unique. Discussions on this topic are
given in \cite{Guttmann1996, conway_hex}. In \cite{conway_hex}, it is
suggested to anisotropise directed site-animals by counting them
according to the number of sites and the number of sites
\emph{supported only from the south} --- \emph{ie} the number of sites
that have a neighbour to the south, but not the east. An example is
given in Figure~\ref{fig dba-aniso eg}. This results \cite{MBM1998} in
the following algebraic generating function
\begin{equation}
  S(q,s) = \frac{1}{2} \left(
    \sqrt{1 - \frac{4q}{(1+q)(1+q-qs)} } - 1
  \right).
\end{equation}

Expanding both $S(q,s)$ and $B(x,y)$ as power series in $s$ and $y$
(respectively) one observes a marked difference in the structure of
their coefficients which leads to our main result. Rewriting $S(q,s) =
\sum_{n \geq 0} R_n(q) s^n$, we find that the $R_n(q)$ are rational
functions of $q$, and the first few are:
\begin{eqnarray}
  R_0(q) & = & \frac{q}{1-q} \nn
  R_1(q) & = & \frac{q^2}{(1-q)^3} \nn
  R_2(q) & = & \frac{q^3(1+q+q^2)}{(1-q)^5(1+q)} \nn
  R_3(q) & = & \frac{q^4(1+2q+4q^2+2q^3+q^4)}{(1-q)^7(1+q)^2}.
\end{eqnarray}
Expansion shows that the numerators are symmetric, positive and
unimodal and that the denominators are given by $D_n(q) = (1-q)^{2n+1}
(1+q)^{n-1}$. Hence the $R_n(q)$ are only singular at $q=\pm 1$. 

Similarly, the generating function $B(x,y)$ can be rewritten as
$B(x,y) = \sum_{n} H_n(x) y^n$, where $H_n(x)$ counts the number of
directed bond-animals with $n$ vertical
bonds according to the number of horizontal bonds they contain. Using
computer enumeration techniques \cite{Jensen2000-dba} one can find%
\footnote{More precisely, the first hundred (or so) terms of the
  expansion of $H_n(x)$ were fitted using Pad\'e approximants.
  Construction of the approximant does not require many series terms,
  and the other terms serve to ``verify'' the conjectured form. We also
  note that in \cite{ADR_haru} it is proved that $H_n(x)$ is rational
  and bounds are given for the degrees of its
  numerator and denominator. } %
the first few $H_n(x)$:
\begin{eqnarray}
  H_0(x) & = & \frac{1}{1-x} \nn
  H_1(x) & = & \frac{1}{(1-x)^3}\nn
  H_2(x) & = & \frac{1+2x+x^2-x^3}{(1-x)^5(1+x)}\nn
  H_3(x) & = & \frac{1+5x+7x^2+x^3-3x^4-2x^5+x^6}{(1-x)^7(1+x)^2}\nn
  H_4(x) & = & \frac{[1,10,33,53,43,3,-25,-20,1,5,2,-1]}{(1-x)^9(1+x)^3(1+x+x^2)},
\end{eqnarray}
where we have written $[a_0,a_1,\dots,a_n]$ in place of
$a_0+a_1y+\dots + a_n y^n$.

We observe that the $H_n(x)$ are simple rational functions whose
denominators are products of cyclotomic polynomials\footnote{We remind
  the reader that the cyclotomic polynomials are the factors of
  $(1-x^n)$, and in particular $(1-x^n) = \prod_{k|n} \Psi_k(x)$,
  where $\Psi_k(x)$ is the the $k^{th}$ cyclotomic polynomial.}. This
structure is quite general and can be proved using the
\emph{haruspicy} techniques described in \cite{ADR_haru}:
\begin{theorem}[from \cite{ADR_haru}]
  \label{thm Hn nature}
  If $B(x,y) = \sum_{n\geq0} H_n(x) y^n$ is the anisotropic generating
  function of directed bond-animals, then 
  \begin{itemize}
  \item $H_n(x)$ is a rational function, 
  \item the degree of the numerator of $H_n(x)$ cannot be greater than the degree of its
    denominator, and
  \item the denominator of $H_n(x)$ is a product of cyclotomic polynomials.
  \end{itemize}
\end{theorem}
If we look a little further we find that the numerators become
increasingly complicated, but the denominators, which we denote
$D_n(x)$, retain a regular structure. Unlike those of directed
site-animals, the denominators of the coefficients of the directed
bond-animal generating function contain higher and higher order
cyclotomic polynomials, and hence have more and more zeros:
\begin{eqnarray}
    D_5(x) & = & (1-x)^{11}(1+x)^4(1+x+x^2)^2\nn
    D_6(x) & = & (1-x)^{13}(1+x)^5(1+x+x^2)^3(1+x^2)\nn
    D_7(x) & = & (1-x)^{15}(1+x)^6(1+x+x^2)^4(1+x^2)^2\nn
    D_8(x) & = & (1-x)^{17}(1+x)^7(1+x+x^2)^5(1+x^2)^3(1+x+x^2+x^3+x^4).
\end{eqnarray}
This dichotomy between the denominators of solved and unsolved models
is observed in many different lattice models and was suggested as the
basis of a numerical test of ``\emph{solvability}'' by Guttmann and
Enting \cite{Guttmann1996, Guttmann2000} --- if one observes an
increasing number of zeros in the denominators of the coefficients of
the anisotropic generating function then the model is probably not
solvable. One can make this notion of solvability more precise by
relating it to differentiably finite functions:
\begin{theorem}[from \cite{MBM_AR_CH}]
  Let $f(x,y) = \sum_{n\geq0} y^n H_n(x)$ be a D-finite series in $y$
  with coefficients $H_n(x)$ that are rational functions of $x$.
  For $n\geq0$ let $S_n$ be the set of poles of $H_n(x)$, and let $S =
  \bigcup_n S_n$. Then $S$ has only a finite number of accumulation
  points.
\label{thm Dfinite poles}
\end{theorem}

Consequently if the set of zeros of the denominators of the
anisotropic generating function has an infinite number of accumulation
points then the anisotropic generating function is not D-finite.
Unfortunately Theorem~\ref{thm Hn nature} does not give sufficiently
detailed information to prove results about the set of singularities
of the coefficients, $H_n(x)$. Ideally, we would like to prove the
exact form of the denominator, which appears to be
\begin{equation}
  \label{eqn denom form}
  D_n(x) = (1-x)^n \prod_{k=1}^{\lfloor n/2 \rfloor+1} \Psi_k(x)^{n-2k+3},
\end{equation}
however this seems to be extremely difficult\footnote{One can
  probably prove that $D_n(x)$ is a factor of the product on the
  right-hand side of this expression using the techniques described
  \cite{ADR_haru} --- proofs of similar results for self-avoiding
  polygons and general bond animals are given in \cite{ADR_haru} and 
  \cite{ADR_sap}.}. %
Instead we prove a weaker result that is still sufficient:
\begin{theorem}
  \label{thm denom 2k-2}
  The denominator of $H_{2k-2}(x)$ contains a factor of $\Psi_k(x)$
  which does not cancel with the numerator, and so $H_{2k-2}$ is
  singular at the zeros of $\Psi_k(x)$.
\end{theorem}
This result them implies:
\begin{cor}
  The singularities of the coefficients $H_n(x)$ in the anisotropic
  generating function $B(x,y)$ form a dense set on the unit circle
  $|x|=1$, and so $B(x,y)$ is not a D-finite power series in $y$.
\end{cor}
Since the specialisation of any D-finite power series is itself
D-finite (provided the specialisation is well-defined --- \emph{ie}
non-singular), we are able to extend this result to directed
bond-animals on any hypercubic lattice.


\section{The proof of Theorem~\ref{thm denom 2k-2}}
\label{sec dfinite proof}
The haruspicy techniques in \cite{ADR_haru} give a way of linking the
``topology'' (in some loose sense) of subsets of bond-animals to the
structure of their generating functions --- and in particular a way of
determining which ``topologies'' cause which singularities. The
following theorem makes this idea precise:
\begin{theorem}[from \cite{ADR_haru}]
  \label{thm poles sections}
  Let $\mathcal{A}_n$ be a dense set of animals with $n$ vertical
  bonds. And let 
  \begin{displaymath}
    H_n(x) = \sum_{A \in \mathcal{A}_n} x^\hhp{A}.
  \end{displaymath}
  If $H_n(x)$ has a denominator factor $\Psi_k(x)$, then there must be
  a section-minimal animal in $\mathcal{A}_n$ that contains a
  $K$-section for some $K \in \mathbb{Z}^+$ divisible by $k$. Further
  if $H_n(x)$ has a denominator factor $\Psi_k(x)^{\alpha}$, then
  there must be a section-minimal animal in $\mathcal{A}_n$ that
  contains $\alpha$ sections that are $K$-sections for some (possibly
  different) $K \in \mathbb{Z}^+$ divisible by $k$.
\end{theorem}
\noindent We have not given definitions of \emph{dense}, \emph{section} and
\emph{section-minimal animal} in the main body of the paper and we
refer the reader to Appendix~\ref{sec haru} (or to \cite{ADR_haru}). Also note
that for convenience we write ``\emph{animal}'' instead of
``\emph{directed bond-animal}''.

\subsection{Animals that cause $\Psi_k(x)$.}

Theorem~\ref{thm denom 2k-2} asserts that a factor of $\Psi_k(x)$
occurs in the denominator of $H_{2k-2}(x)$. According to the above
theorem this can only be the case if there is a section-minimal animal
with $2k-2$ vertical bonds that contains at least one $k$-section (or
a $K$-section with $K$ an integer multiple of $k$). We start by
characterising such animals.

\begin{figure}[ht]
  \centering
  \includegraphics[scale=0.6]{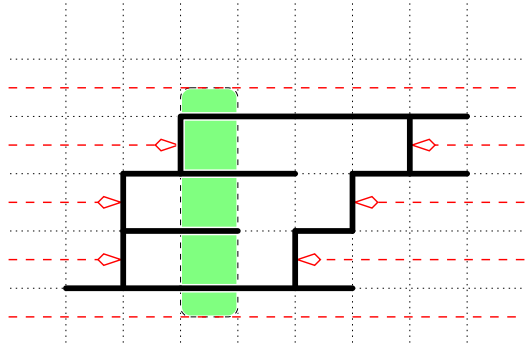}
  \caption{A directed bond-animal with $2k-2$ vertical bonds
    and a $k$-section (highlighted).}
  \label{fig dba ksec}
\end{figure}

\begin{lemma}
  Let $A$ be an animal that contains a $k$-section. $A$ must contain
  at least $2k-2$ vertical bonds. If $A$ contains a $k$-section and
  exactly $2k-2$ vertical bonds then there must be exactly $2$
  vertical bonds in each row of $A$.
  \label{lem first ksec dba}
\end{lemma}
\proof Consider an animal that contains a $k$-section. The $k$-section
must contain at least $k-1$ cells in a vertical line (see
Figure~\ref{fig dba ksec}). In order to be a $k$-section, no
section-line may cross any of these cells. Hence each section line to
the left and right of these cells must be obstructed by a vertical
bond and so there must be at least $1$ vertical bond to the left and
$1$ vertical bond to the right of each of these cells. Hence an animal
that contains a $k$-section must contain at least $2k-2$ vertical bonds.

By similar reasoning, if the animal contains exactly $2k-2$ vertical
bonds then there must be 2 vertical bonds in each row. \qed

We note that one can push the above proof further to show that a
directed bond-animal with exactly $2k-2$ vertical bonds contains no
more than one $k$-section, however we do not need this result. We also
note that the above lemma and Theorem~\ref{thm poles sections} imply
that the denominators of $H_n(x)$ with $n<2k-2$ cannot contain a
factor of $\Psi_k(x)$.

\spacebreak

The previous Lemma shows that the factor of $\Psi_k(x)$ in the
denominator of $H_{2k-2}(x)$ is caused by those section-minimal
animals that contain a $k$-section, which requires that they have $2$
vertical bonds in each row. In order to prove that this denominator
factor does not cancel with the numerator of $H_{2k-2}(x)$, we need to
examine the set of all directed bond-animals with $2$ vertical bonds
per row.

\begin{figure}[ht]
  \centering
  \includegraphics[scale=0.5]{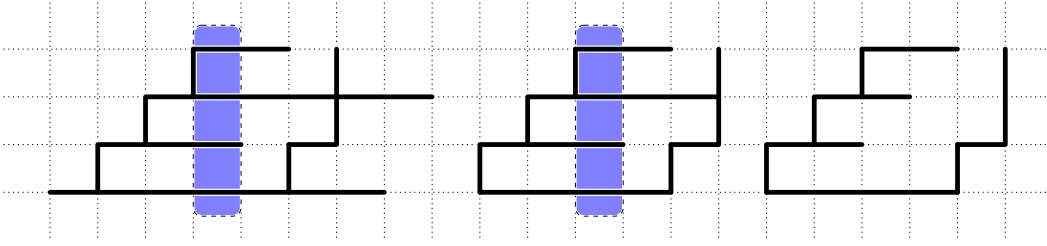}
  \caption{(left) A \emph{2}-directed bond-animal which has $2$ vertical bonds in each
    row. (centre) The corresponding \emph{primitive 2}-directed bond
    animal. (right) A \emph{2}-directed animal that contains $2k-2$
    vertical bonds but no $k$-section.}
  \label{fig 2directed eg}
\end{figure}

\begin{defn}
  A \emph{2-directed animal} is a directed animal which has $2$
  vertical bonds in each row. A \emph{primitive 2-directed animal} is
  a \emph{2}-directed animal in which all vertices of degree~$1$ lie
  between vertical bonds (see Figure~\ref{fig 2directed eg}).
\end{defn}
By Lemma~\ref{lem first ksec dba} all animals that contain a
$k$-section and have $2k-2$ vertical bonds are \emph{2}-directed
animals, but there are \emph{2}-directed animals with $2k-2$ vertical
bonds that do not contain a $k$-section (see Figure~\ref{fig 2directed
  eg}).  Also, one may construct a \emph{2}-directed animal from a
\emph{primitive 2}-directed animal by prepending a line of horizontal
bonds to the left of the bottom-leftmost-vertex, and appending lines
of horizontal bonds to the right of the rightmost vertices.
Consequently, if $f_n(x)$ is the generating function of
\emph{primitive 2}-directed animals with $2n$ vertical bonds, then
$\left( \frac{1}{1-x} \right)^{n+2} f_n(x)$ is the generating function
of all \emph{2}-directed animals with $2n$ vertical bonds.

\begin{lemma}
  \label{lem 2dba ok}
  The generating function of \emph{2}-directed animals with $2n-2$
  vertical bonds ($n>0$) has poles at the zeros of $\Psi_{n}(x)$ if
  and only if $H_{2n-2}(x)$ has poles at the zeros of $\Psi_{n}(x)$.
\end{lemma}
\proof Since section deletion and duplication do not alter the number
of vertical bonds, nor move them between rows, it follows that
\emph{2}-directed animals are closed under section duplication and
deletion and so form a dense set. Similarly the set of directed bond
animals that are not \emph{2}-directed bond-animals is dense. This
means that we may apply Theorem~\ref{thm poles sections} to both of these sets.

Let $\mathcal{A}$ be the set of all directed bond-animals with $2k-2$
vertical bonds, and let $\mathcal{B}$ be the set of all
\emph{2}-directed animals with $2k-2$ vertical bonds. Now split
$H_n(x)$ into a sum over the animals in $\mathcal{B}$ and all the others:
\begin{eqnarray}
  H_{2n-2}(x) & = & \sum_{A \in \mathcal{B}} x^\hhp{A} 
  + \sum_{A \in \mathcal{A} \setminus \mathcal{B}} x^\hhp{A} \\
  & = & G_1(x) + G_2(x). \nonumber
\end{eqnarray}
By Theorem~\ref{thm Hn nature} we know that $G_1(x)$ and $G_2(x)$ are
rational generating functions whose denominators are products of
cyclotomic polynomials. Since all those section-minimal animals with
$k$-sections contribute to $G_1(x)$ and not $G_2(x)$, by
Theorem~\ref{thm poles sections} there is no factor of $\Psi_k(x)$ (or
higher cyclotomic factors) in the denominator of $G_2(x)$.

Let $G_1(x)$ have a factor of $\Psi_k(x)^{\alpha}$ in its denominator
that does not cancel with its numerator.  Since there are no factors
of $\Psi_k(x)$ in the denominator of $G_2(x)$, it follows that
$H_{2k-2}(x)$ also has a factor of $\Psi_k(x)^\alpha$ in its
denominator. Similarly if $H_n(x)$ has a factor of $\Psi_k(x)^\alpha$
in its denominator that does not cancel with its numerator, then so
must $G_1(x)$. \qed

The above lemma makes the proof of Theorem~\ref{thm denom 2k-2} much
simpler. Instead of having to analyse all directed bond-animals, we
only need look at a much simpler subset --- \emph{2}-directed
animals. Further we don't have to enumerate this subset exactly, we
only need to locate the singularities of its generating function. 

\subsection{Counting \emph{2}-directed animals}

In order to study the generating function of \emph{2}-directed animals
we make use of a powerful enumeration technique, the Temperley method.
The method consists (essentially) of two steps --- finding a
recurrence satisfied by coefficients or generating functions, and then
solving that recurrence. For the purposes of this paper we need
to analyse the singularities of the generating function, and it
transpires that an expression for the generating function is
unnecessary --- it is sufficient to work with the recurrences it
satisfies. As was the case in \cite{ADR_sap} we use a
variation of the Temperley method involving Hadamard products. 

We start by defining the restricted Hadamard product and then showing
how it may be used to find a recurrence satisfied by the generating
function of \emph{2}-directed animals.
\begin{defn}
  Let $f(t) = \sum_{n \geq 0} f_n t^n$ and $g(t)=\sum_{n \geq 0} g_n
  t^n$ be formal power series in $t$. The (restricted) Hadamard
  product is defined to be
  \begin{displaymath}
    f(t) \hp{t} g(t) = \sum_{n \geq 0} f_n g_n.
  \end{displaymath}
  We note that if  $f(t)$ and $g(t)$ are two power series with real
  coefficients such that 
  \begin{displaymath}
    \lim_{n \to \infty} \left| f_n g_n \right|^{1/n} < 1,    
  \end{displaymath}
  then the Hadamard product $f(t) \hp{t} g(t)$ will exist.
\end{defn}
Below we consider Hadamard products of power series in $t$ whose
coefficients are power series in two variables $x$ and $s$. The
products are of the form $f(t;x) \hp{t} T(t,s;x) = \sum_{n \geq 0}
f_n(x) T_n(s;x)$. The summands are the generating functions of certain
directed bond animals and it follows that the $n^{th}$ summand is $O(sx^n)$
and so the sum converges. In order to re-express the Hadamard products
we will use the following result:
\begin{lemma}
  \label{lem easy hadamard}
  Let $f(t)$ be a formal power series in $t$. The following
  (restricted) Hadamard products are easily evaluated.
  \begin{eqnarray*}
    f(t) \hp{t} \frac{1}{1-\alpha t} & = & f(\alpha) \\
    f(t) \hp{t} \frac{ n! t^n }{(1-\alpha t)^{n+1}} & = &
    \left.\Ddiff{f}{t}{n} \right|_{t = \alpha}
  \end{eqnarray*}
  We also note that the Hadamard product is a linear operator.
\end{lemma}
\proof See similar lemma in \cite{ADR_sap}.

\spacebreak

Every \emph{2}-directed animal may be constructed row by row --- many
other objects have been counted in this way. In this paper we use the
same variation of this technique used in \cite{ADR_sap} which involves
decomposing the object into a \emph{seed} and \emph{building blocks}.
To simplify the following discussion we will work with primitive
\emph{2}-directed animals rather than all \emph{2}-directed animals;
since their generating functions differ only by factors of $(1-x)$,
the other cyclotomic factors are unaffected. For convenience we shall
drop the word ``\emph{primitive}''.

\begin{figure}[ht]
  \centering
  \includegraphics[scale=0.4]{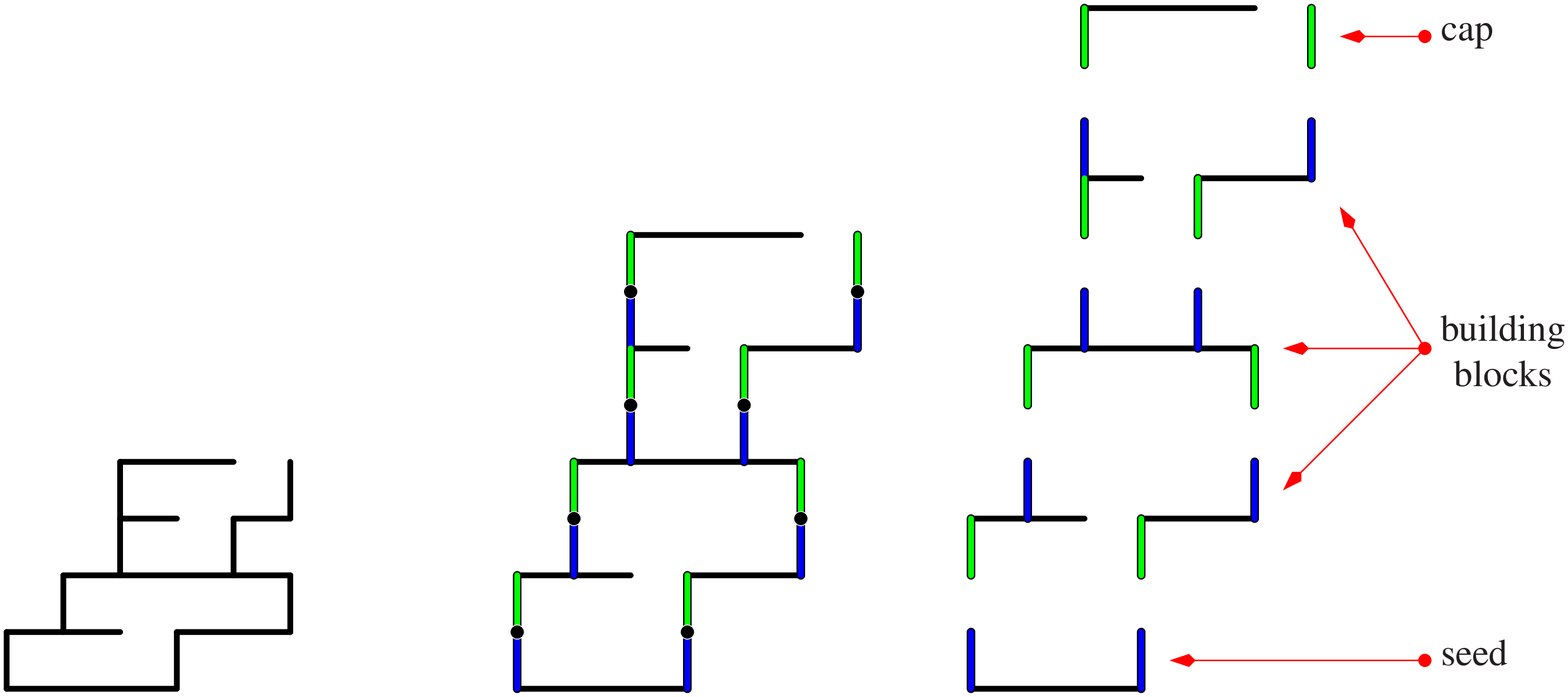}
  \caption{Decomposing a \emph{2}-directed animal into a seed, a
    sequence of building blocks and then a cap.}
  \label{fig dba decomp}
\end{figure}

Start with a \emph{2}-directed animal and duplicate every row
(including the vertical bonds in each row) --- see Figure~\ref{fig dba
  decomp} (left and centre). Now cut horizontally through the centre
of each pair of duplicated rows; this decomposes the animal into a
``\emph{seed block}'' (occupying a single row at the bottom of the
animal), a sequence of ``\emph{building blocks}'' (each occupying two
rows) and then a ``\emph{cap}'' (occupying a single row at the top of
the animal) --- see Figure~\ref{fig dba decomp} (top). We note that
the sequence of blocks is restricted so that the top row of one block
must have the same length as the bottom row of the next block --- the
Hadamard product allows us to easily translate this restriction into
an operation on generating functions.

We are able to find a recurrence satisfied by the generating function
of \emph{2}-directed animals from the generating functions of the
seeds, building blocks, and caps. In particular we must enumerate each
of these objects according to the number of horizontal bonds, and the
distance between the vertical bonds.

The seed is simply a line of horizontal bonds terminated on each end
by a vertical bond. It has generating function $\frac{sx}{1-sx}$
(where $s$ is conjugate to the distance between the vertical bonds).

\begin{figure}[ht]
  \begin{center}
  \includegraphics[scale=0.5]{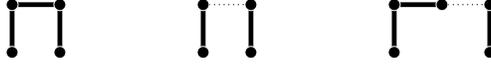}
  \caption{The section-minimal caps.}
  \label{fig dba cap}
\end{center}
\end{figure}

The caps consist of two vertical bonds with some number of horizontal
bonds between them. Since the animal is directed, these horizontal
bonds must be attached to the left-hand vertical bond, but not
necessarily the right-hand vertical bond. The section-minimal caps are
given in Figure~\ref{fig dba cap}, and expanding them gives the
generating function:
\begin{equation}
  \label{eqn kap gf}
  \frac{t(1+x-tx)}{(1-t)(1-tx)} = -1 + \frac{1}{(1-x)(1-t)} - \frac{x}{(1-x)(1-tx)},
\end{equation}
where $t$ is conjugate to the distance between the vertical bonds.

\begin{figure}[ht]
  \centering
  \psfrag{where}{where} \psfrag{is}{is} \psfrag{or}{or}
  \psfrag{sothat}{so that} \psfrag{and}{and}
  \includegraphics[scale=0.5]{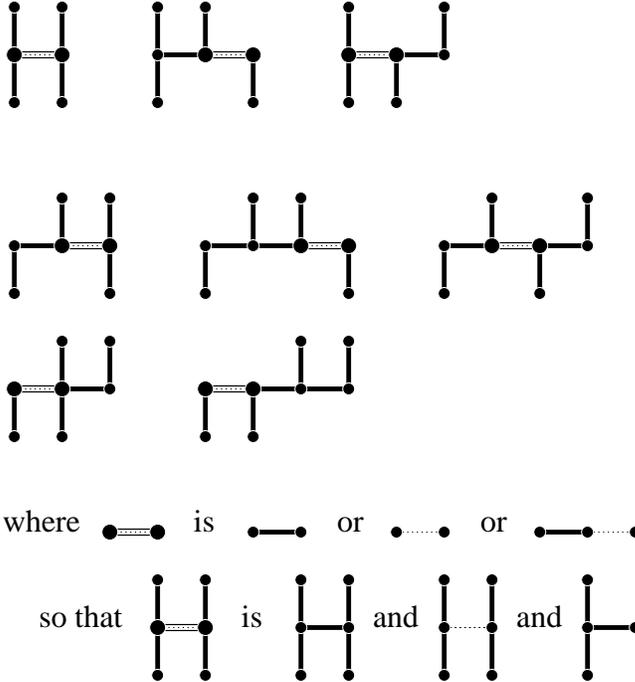}
  \caption{The section-minimal building blocks of \emph{2}-directed
    bond-animals. The highlighted horizontal bonds are short-hand for
    either a single horizontal bond, a single horizontal space (with
    no bond) or a single horizontal bond and a single horizontal space.}
  \label{fig dba basic}
\end{figure}

The building blocks are (reasonably) complicated and we give the
section-minimal building blocks in Figure~\ref{fig dba basic}. We
compute the generating function of the building blocks by expanding
each of the sections. We now need two extra variables: $s$ and $t$ are
conjugate to the distances between the vertical bonds in the top and
bottom rows (respectively). Doing this gives (moving from
left-to-right and top-to-bottom in Figure~\ref{fig dba basic}):
\begin{eqnarray}
  T(s,t;x)  & = & \Big(\seq{stx} + \seq{st} + \seq{stx} \seq{st}\Big) 
  + \seq{stx}\Big(\seq{tx} + \seq{t} + \seq{tx}\seq{t} \Big) \nn
  && + \Big(\seq{stx} + \seq{st} + \seq{stx}\seq{st} \Big) \seq{sx} \nn
  && + \seq{tx} \Big( \seq{stx} + \seq{st} + \seq{stx}\seq{st} \Big)
  + \seq{tx} \seq{stx} \Big(\seq{tx} + \seq{t} + \seq{tx}\seq{t}\Big)\nn
  && + \seq{tx} \Big(\seq{stx} + \seq{st} + \seq{stx}\seq{st} \Big) \seq{sx} \nn
  && + \Big(\seq{tx} + \seq{t} + \seq{tx}\seq{t} \Big)\seq{sx}
   + \Big(\seq{tx} + \seq{t} + \seq{tx}\seq{t} \Big)\seq{x}\seq{sx}
\end{eqnarray}
where we have used the short hand notation $\seq{f} = \frac{f}{1-f}$.
This may then be written in (a slightly non-standard) partial fraction
form as:
\begin{multline}
  \label{eqn bb gf1}
  T(s,t;x)  = - \frac{sx}{(1-x)(1-sx)} \Big( t^0 \Big) 
  + \frac{2sx}{(1-x)^2(1-sx)} \left(\frac{1}{1-t} \right) \\
  + \frac{s}{(1-x)(1-sx)(s-x)} \left(\frac{1}{1-st} \right)\\
  +  \frac{(s-1)s - (s-2)(s^2-s+1)sx - (s^2-s+3)s^2x^2+(2s^2+1)sx^3 -s^2x^4}
{(1-x)^2(1-s)^2(1-sx)(s-x)}\left(\frac{1}{1-tx} \right) \\
  - \frac{sx^2}{(1-x)(1-s)} \left(\frac{t}{(1-tx)^2} \right) + \frac{s(1 - (1+x-x^2)s)}{(1-x)(1-s)^2(1-sx)} \left(\frac{1}{1-stx} \right),
\end{multline}
which we shall rewrite (more concisely) as:
\begin{multline}
  \label{eqn bb gf2}
  T(s,t;x)  = c_0 \Big( t^0 \Big) + c_1 \Big(\frac{1}{1-t} \Big) 
  + c_2  \Big(\frac{1}{1-st} \Big) + c_3 \Big(\frac{1}{1-tx} \Big)
  + c_4 \Big(\frac{t}{(1-tx)^2} \Big) + c_5 \Big(\frac{1}{1-stx} \Big).
\end{multline}
where the $c_i$ are the corresponding rational functions of $s$ and
$x$.

We find a recurrence for \emph{2}-directed animals in two
steps. First we consider ``\emph{uncapped}'' \emph{2}-directed
animals, which are those generated from a seed and sequence of
building blocks, but no cap. These are simply \emph{2}-directed
animals with no horizontal bonds attached to the top of vertical bonds
in their topmost row. We then find a recurrence for all
\emph{2}-directed animals in terms of the uncapped \emph{2}-directed animals.

\begin{lemma}
  \label{lem uncapped rec}
  Let $\breve{f}_n(s;x)$ be the generating generating function of
  \emph{uncapped} \emph{2}-directed bond-animals (with $2n$ vertical
  bonds). The variable $x$ is conjugate to the number of horizontal
  bonds and $s$ is conjugate to the distance between the vertical
  bonds in the topmost row. This generating function satisfies the
  following functional equation:
  \begin{eqnarray}
    \breve{f}_1(s;x) & = & \frac{sx}{1-sx} \\
    \breve{f}_{n+1}(s;x) & = &  c_1 \breve{f}_n(1;x)
    + c_2 \breve{f}_n(s;x) +c_3 \breve{f}_n(x;x)
    + c_4 \left. \Ddiff{ \breve{f}_n}{s}{} \right|_{s=x} 
    + c_5 \breve{f}_n(sx;x),
  \end{eqnarray}
  where the $c_i$ are given in equations~\Ref{eqn bb gf1} and~\Ref{eqn bb gf2}.
  
  This recurrence is singular at two points of interest, namely $s=1$
  and $s=x$. At these points singularities of the building block
  generating function coalesce and the recurrences change structure:
  \begin{eqnarray}
    \breve{f}_{n+1}(1;x) & = &  \frac{1+2x}{(1-x)^2} \breve{f}_n(1;x)
    - \frac{1+x+x^2}{(1-x)^2} \breve{f}_n(x;x) 
    - \frac{x(1+x)}{(1-x)^2} \left. \Ddiff{ \breve{f}_n}{s}{} \right|_{s=x} \nn
    && - \frac{x^3}{2(1-x)} \left. \Ddiff{ \breve{f}_n}{s}{2} \right|_{s=x}\\ 
    \breve{f}_{n+1}(x;x) & = &  \frac{2x^2}{(1-x)^2(1-x^2)}  \breve{f}_n(1;x)
    - \frac{x}{(1-x)^3}\breve{f}_n(x;x) 
    - \frac{x(1-x^2-x^3)}{(1-x)(1-x^2)} \left. \Ddiff{ \breve{f}_n}{s}{} \right|_{s=x} \nn
    && + \frac{x}{(1-x)^2} \breve{f}_n(x^2;x).
  \end{eqnarray}
\end{lemma}
\proof The generating function of uncapped \emph{2}-directed animals
with $2$ vertical bonds is exactly that of the seed generating
function, namely $\frac{sx}{1-sx}$. We then obtain the generating
functions, $\breve{f}_n(s;x)$, by repeatedly adding building
blocks. 

Let $\breve{f}_n(s;x) = \sum_{m \geq 1} \breve{f}_{n,m}(x) s^m$, and
$T(s,t;x) = \sum_{m \geq 1} T_m(s;x) t^n$. The coefficient
$\breve{f}_{n,m}(x)$ counts those \emph{2}-directed animals which have
$m$ cells separating the $2$ vertical bonds in their \emph{top} row.
Similarly $T_m(s;x)$ counts those building blocks with $m$ cells
separating the $2$ vertical bonds in their \emph{bottom} row. Thus
adding a new building block corresponds to the following operation on
the generating functions:
\begin{eqnarray*}
  \breve{f}_{n+1}(s;x) & = & \sum_{m \geq 1} \breve{f}_{n,m}(x) T_m(s;x) \\
  & = & \breve{f}_n(t;x) \hp{t} T(s,t;x).
\end{eqnarray*}
Applying Lemma~\ref{lem easy hadamard} to the partial fraction form of
$T(s,t;x)$ gives the first recurrence. Repeating this with $s=1$ and
$s=x$ gives the later recurrences. Note that $\breve{f}_n(0;x) = 0$,
since there must be some positive number of cells separating the
vertical bonds in the top row of the animal.\qed

\begin{lemma}
  \label{lem kap rec}
  The generating function, $f_n(x)$, of \emph{all} \emph{2}-directed
  animals with $2n$ vertical bonds may be expressed in terms of the
  generating function of uncapped \emph{2}-directed animals:
  \begin{equation}
    f_{n}(x)  = \frac{1}{1-x} \big(\breve{f}_n(1;x) - x \breve{f}_n(x;x) \big).
  \end{equation}
\end{lemma}
\proof By similar reasoning to that given in the proof of the previous
lemma, we may express the capped generating function as a Hadamard
product of the uncapped generating functions together with the
generating function of the caps (see equation~\Ref{eqn kap gf}). Again
we make use of the fact that $\breve{f}_n(0;x) = 0$. \qed

\subsection{Analysing the singularities}
Using the recurrences for the generating functions of
\emph{2}-directed bond-animals we proceed in two steps. We iterate the
recurrences in order to determine the structure of the coefficients as
rational functions of $s$ and $x$. We then substitute this structure
back into the recurrence to link the singularities of the generating
function of animals with $2n$ vertical bonds at $s=1$ to those with
$2n-2$ vertical bonds at $s=x$. Continuing this reasoning, we link the
singularities of $f_n$ --- a function we do not know in closed form
--- to the singularities of $f_1$ --- which is a simple rational
function that we do know.

Examining the first few generating functions, $\breve{f}_n(s;x)$ we
see that their denominators are products of cyclotomic
polynomials, $\Psi_k(x)$, and factors of the form $(1-sx^k)$. To refer
easily to polynomials of this type we define the following sets:
\begin{defn}
  Let $\mathbb{C}_n(s;x)$ be the set of all polynomials of the form
\begin{equation}
  \prod_{k=1}^n (1-sx^k)^{a_k} \Psi_k(x)^{b_k},
\end{equation}
where $a_k$ and $b_k$ are non-negative integers. We also define
$\mathbb{C}_n(x) = \mathbb{C}_n(0;x)$.
\end{defn}
Using the above notation we can describe the structure of these
generating functions:
\begin{lemma}
  The generating function, $\breve{f}_n(s;x)$, is of the form:
  \begin{equation}
    \breve{f}_n(s;x) = \frac{N_n(s;x)}{D_n(s;x) (1-sx^{n})}
  \end{equation}
  where $N_n(s;x)$ and $D_n(s;x)$ are polynomials in $s$ and $x$,
  with the further restriction that $D_n(s;x) \in \mathbb{C}_{n-1}(s;x)$.
  \label{lem uncapped form}
\end{lemma}
\proof We first note that since $\breve{f}_n(s;x)$ counts uncapped
\emph{2}-directed animals with $2n$ vertical bonds, it cannot be
singular at $s=1$, and so its denominator does not contain factors of
$(1-s)$. The result then follows by iteration of the recurrence. See
\cite{ADR_sap, MBM_ADR_02} for similar arguments. \qed

Before we can substitute the above form into the recurrences satisfied
by $\breve{f}_n$, we need to show that one of the coefficients of the
recurrence does not have zeros on the unit circle which could
potentially cancel singularities of $\breve{f}_n$.
\begin{lemma}
  \label{lem c5 zero}
  At $s = x^n$, the coefficient $c_5(x^n;x) =
  \frac{x^n(1-(1+x-x^2)x^n)}{(1-x)(1-x^n)^2(1-x^{n+1})}$ is non-zero
  everywhere on the unit-circle $|x| = 1$.
\end{lemma}
\proof Consider the zeros of the numerator polynomial $(1 -
(1+x-x^2)x^n)=0$. We may rewrite this as
\begin{equation}
  x^n = \frac{1}{1+x-x^2}.
\end{equation}
If the polynomial has a zero on the unit circle, $x = e^{i \theta}$,
then it follows that $| 1+x-x^2 | =1$. This then gives
\begin{equation}
  \big(1+ \cos(\theta) + \cos(2\theta) \big)^2 + \big(\sin(\theta) + \sin(2\theta) \big)^2
  = 1,
\end{equation}
which reduces to the condition $\cos^2(\theta) = 1$. Hence the only
candidates for zeros are $x = \pm 1$. Inspection of the polynomial
then shows that it has a single zero at $x=1$ for all $n$, and that it
has a single zero at $x=-1$ for odd $n$.

Since the denominator of $c_5(x^n,x)$ contains factors of $(1-x)$ and
$(1+x)$ for all $n \geq 1$, it follows that neither $x=1$ or $x=-1$ is
a zero of the function. \qed

\vspace{0.5cm}

\begin{theorem}
  \label{thm dba sing}
  For all $n \geq 1$, the generating function $f_n(x)$ has simple
  poles at the zeros of $\Psi_{n+1}(x)$.
\end{theorem}
\proof Fix $n$ and let $\xi$ be a zero of $\Psi_{n+1}(x)$. We will
start by showing that $\breve{f}_{k}(x^{n-k+1};x)$ is singular at
$x=\xi$ by induction on $k$ for fixed $n$.  We then show that this is
sufficient to prove the above theorem by linking the singularities of
$f_n$ to those of $\breve{f}_n$.

\vspace{0.2cm}

\noindent Setting $k=1$ gives $\breve{f}_1(x^n;x) = \frac{x^{n+1}}{1-x^{n+1}}$.
which is singular at $x = \xi$.
  
\vspace{0.2cm}

\noindent We now proceed by induction on $k$ using the recurrences of
Lemma~\ref{lem uncapped rec}. Assume that $\breve{f}_{k}(x^{n-k+1};x)$
is singular at $x=\xi$.  By Lemma~\ref{lem uncapped form}, we may
write $\breve{f}_n(s;x)$ as:
\begin{displaymath}
  \breve{f}_k(s;x) = \frac{N_k(s;x)}{D_k(s;x)(1-sx^{k})},
\end{displaymath}
where $D_k(s;x) \in \mathbb{C}_{k-1}(s;x)$ and $N_k(s;x)$ is some
polynomial in $s$ and $x$. Substitute this form into the recurrences
of Lemma~\ref{lem uncapped rec}. We may now write
$\breve{f}_{k+1}(s;x)$ may as
\begin{equation}
  \breve{f}_{k+1}(s;x) = \frac{N(s;x)}{D(s;x)} + c_5(s;x) \breve{f}_k(sx;x),
\end{equation}
where $N(s;x)$ and $D(s;x)$ are polynomials and $D(s;x) \in
\mathbb{C}_{k}(s;x)$. Setting $s=x^{n-k}$ gives
\begin{equation}
  \breve{f}_{k+1}(x^{n-k};x) = \frac{N(x^{n-k};x)}{D(x^{n-k};x)} 
  + c_5(x^{n-k};x) \breve{f}_k(x^{n-k+1};x),
\end{equation}
with $D(x^{n-k};x) \in \mathbb{C}_{n}(x)$. By Lemma~\ref{lem c5 zero} we know
that $c_5(x^{n-k};x)$ is not zero at $x=\xi$. Since
$\breve{f}_k(x^{n-k+1};x)$ is singular at $x=\xi$ so is
$\breve{f}_{k+1}(x^{n-k};x)$. By induction we have shown that
$\breve{f}_n(x;x)$ is singular at $x=\xi$. Further, Lemma~\ref{lem
  uncapped form} implies that the singularity is a simple pole.

\vspace{0.3cm}

\noindent Using Lemma~\ref{lem kap rec} the singularities of $\breve{f}_n$ are
linked to those of $f_n$:
\begin{equation}
  f_{n}(x) = \frac{1}{1-x} \left( \breve{f}_n(1;x) - x\breve{f}_n(x;x) \right).
\end{equation}
Lemma~\ref{lem uncapped form} then implies that $\breve{f}_n(1;x)$ is
not singular at $x = \xi$ and so the simple pole of $\breve{f}_n(x;x)$
at $x = \xi$ implies a simple pole in $f_n(x)$ at $x = \xi$. \qed

\spacebreak

The above theorem gives our main result:
\begin{cor}
  Since $f_n(x)$ has simple poles at the zeros of $\Psi_{n+1}(x)$,
  the coefficient $H_{2n}(x)$ has simple poles at the zeros of
  $\Psi_{n+1}(x)$ and the anisotropic generating function of directed
  bond-animals is not differentiably finite.
\end{cor}
\proof Since the generating function of primitive \emph{2}-directed
animals and \emph{2}-directed animals are related by factors of
$(1-x)$, it follows from Theorem~\ref{thm dba sing} and Lemma~\ref{lem
  2dba ok} that $H_{2n}(x)$ has simple poles at the zeros of
$\Psi_{n+1}(x)$.

Let $S$ be the union of the singularities of $H_n(x)$ for all $n$. For
any $q \in \mathbb{Q}$ there exists $k$ such that $\Psi_k(e^{2\pi i
  q}) = 0$, and since $H_{2k-2}(x)$ has simple poles at the zeros of
$\Psi_k(x)$, it follows that $e^{2 \pi i q} \in S$. Consequently $S$
is dense on the unit circle $|x|=1$, and by Theorem~\ref{thm Dfinite
  poles} the anisotropic generating function of directed bond-animals
is not differentiably finite. \qed

This can then be extended to give the following result:
\begin{cor}
  \label{cor dfinite hypercube}
  Let $\mathcal{B}_d$ be the set of directed bond-animals on the
  $d$-dimensional hypercubic lattice, and let $B_d$ be the anisotropic
  generating function
  \begin{displaymath}
    B_d(x_1, \dots, x_{d-1}, y) = 
    \sum_{A \in \mathcal{B}_d} y^{|A|_d} \prod_{i=1}^{d-1} x_i^{|A|_i},
  \end{displaymath}
  where $|A|_i$ is the number of bonds parallel to the unit vector
  $\vec{e_i}$. Then $B_1(y) = \frac{1}{1-y}$, and for all $d \geq 2$ the
  generating function is not a D-finite power series in $y$.
\end{cor}
\proof When $d=1$ the only directed bond-animals consist of a line of
bonds; the generating function is simply $\frac{1}{1-y}$. When $d=2$ the
result follows from the previous corollary. Finally if $d>2$, set $x_2
= \dots = x_{d-1} = 0$ in $B_d$. This specialisation is well-defined
since the generating function now counts those animals that are
confined to the plane spanned by $\{ \vec{e_1}, \vec{e_d} \}$ which
are simply directed bond-animals on the square lattice.

Since the well defined specialisation of a D-finite power series is
itself D-finite \cite{Lipshitz1989}, it follows that if $B_d$ were a
D-finite function of $y$, then $B(x,y)$ would also be D-finite.  This
contradicts the previous corollary and the result follows. \qed


\section{Conclusion}
We have demonstrated that the anisotropic generating function of
directed bond-animals is not differentiably finite and so is
fundamentally different from that of directed site-animals which has
been solved. 

Unfortunately this result does not enable us to say anything rigorous
about the nature of the isotropic generating function; one can readily
construct an example of a function, $f(x,y)$ which is not D-finite
that becomes D-finite when $x=y$. For example:
\begin{equation}
  F(x,y) = \sum_{n \geq 1} \frac{y^n}{(1-x^n)(1-x^{n+1})}.
\end{equation}
is not a D-finite function of $y$ by Theorem~\ref{thm Dfinite
  poles}. Setting $y=x=z$ gives a rational, and hence D-finite,
function of $z$:
\begin{equation}
  F(z,z) = \sum_{n \geq 1} \frac{z^n}{(1-z^n)(1-z^{n+1})} 
  = \frac{z}{(1-z)^2}.
\end{equation}

On the other hand, the ``\emph{anisotropisation}'' of models that have
been solved does alter the nature of their generating functions,
rather it moves singularities around in the complex plane. Of course,
this does not imply anything about unsolved problems.  It should also
be noted that there exist non-rigorous Renormalisation Group arguments
which imply that anisotropisation should not affect the analytic
nature of the generating function \cite{renorm}. We note that if the
isotropic generating function is indeed not D-finite then it will not
be found using computer packages such as GFUN \cite{Algolib} or
differential approximants \cite{Guttmann1989} which can only find
D-finite solutions.

We are currently working on extending non-D-finiteness results to
other bond-animal problems including square lattice bond-animals and
bond-trees.  Unfortunately, work on a similar result for self-avoiding
walks appears to be beyond the scope of these techniques \cite{ANR}
--- the self-avoiding walk analogue of \emph{2}-directed animals and
\emph{2-4-2} polygons (see \cite{ADR_haru}) appear to be quite
complicated and so finding recurrences such as those in Lemma~\ref{lem
  uncapped rec} would be very difficult.

Finally, it may also be possible to extend the haruspicy techniques to
site-animals and polyominoes making it possible to show that
self-avoiding polygons or general site-animals, counted by an
``\emph{anisotropised}'' area are not D-finite. This would also
possibly explain why directed site-animals on the hexagonal lattice
remain unsolved --- there is strong numerical evidence
\cite{conway_hex} indicating that their anisotropic generating
function is not D-finite, and it may be possible to sharpen this
evidence into proof.

\section*{Acknowledgements}
I would like to thank I. Jensen for his anisotropic series data and A.
J. Guttmann for his help with the manuscript.  This work was partially
funded by the Australian Research Council.

\appendix
\section{Haruspicy}
\label{sec haru}

In a previous paper \cite{ADR_haru} haruspicy\footnote{The word
  ``haruspicy'' refers to techniques of divination based on the
  examination of the forms and shapes of the organs of animals.}
techniques have been developed that allow us to determine properties
of the anisotropic generating function of a set of bond animals
without detailed knowledge of those animals. This allows the
techniques to applied to problems that are unsolved, such as
self-avoiding polygons \cite{ADR_sap} and (in this paper) directed
bond-animals.

The basic idea is to reduce or squash the set of animals down onto
some minimal set, and then determine properties of the coefficients,
$H_n(x)$ of the anisotopic generating function by examining the bond
configurations of the minimal animals.

\begin{figure}[h!]
  \begin{center}
    \includegraphics[scale=0.5]{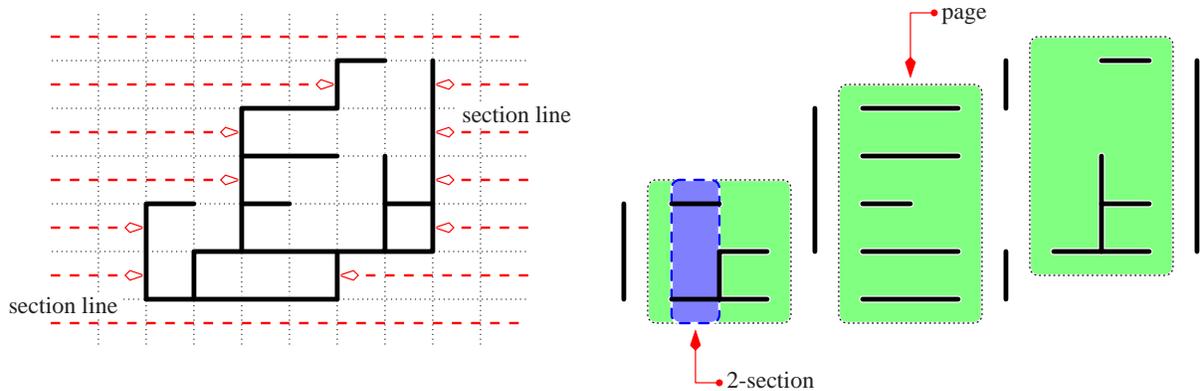}\\
    \caption{\emph{Section lines} (the heavy dashed lines in the
      left-hand figure) split the animal into \emph{pages} (as shown
      on the right-hand figure).  Each column in a page is a
      \emph{section}. This animal is split into 3 pages, each
      containing two sections; a $2$-section is highlighted.  11
      vertical bonds lie between pages and 3 vertical bonds lie within
      the pages.}
    \label{fig section defn}
  \end{center}
\end{figure}

We start by showing how directed bond animals may be cut up so that
they may be ``squashed'' in a consistent way.
\begin{defn}
  Draw horizontal lines from the extreme left and the extreme right of
  the lattice towards the animal so that the lines run through the
  middle of each lattice cell.  These lines are called \emph{section
    lines}.  The lines are terminated when they first touch (\emph{ie}
  are obstructed by) a vertical bond (see Figure~\ref{fig section
    defn}).
  
  Cut the lattice along each section line from infinity until it
  terminates at a vertical bond. Then from this vertical bond cut
  vertically in both directions until another section line is reached.
  In this way the animal (and the lattice) is split into \emph{pages}
  (see Figure~\ref{fig section defn}); we consider the vertical bonds
  along these vertical cuts to lie \emph{between} pages, while the
  other vertical bonds lie \emph{within} the pages.
  
  We call a \emph{section} the set of horizontal bonds within a single
  column of a given page.  Equivalently, it is the set of horizontal
  bonds of a column of an animal between two neighbouring section
  lines. A section with $k$ horizontal bonds is a $k$-section.  The
  number of $k$-sections in an animal, $P$, is denoted by
  $\sigma_k(P)$.
\end{defn}

\begin{figure}[ht]
  \begin{center}
    \includegraphics[scale=0.53]{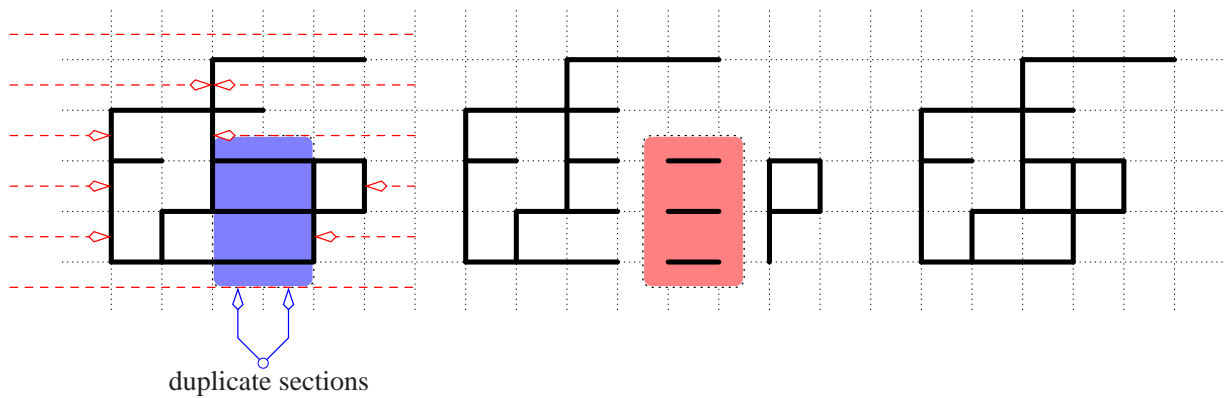}
    \caption{The process of section deletion. The two indicated sections 
      are identical. Slice either side of the duplicate and separate the animal into
      three pieces. The middle piece, being the duplicate, is removed and the remainder of
      the animal is recombined. Reversing the steps leads to section
      duplication.}
    \label{fig explicit section delete}
  \end{center}
\end{figure}

By dividing an animal into sections we see that many of the sections
are superfluous and are not needed to encode its ``\emph{shape}'' (in
some loose sense of the word). In particular, if there are two
identical sections next to each other, then we can reduce the animal
by removing one of them.
\begin{defn}
  We say that a section is a \emph{duplicate section} if the section
  immediately on its left is identical
  (see Figure~\ref{fig explicit section delete}).
  
  An animal can be reduced by \emph{deleting} duplicate sections;
  slice the animal on either side of the duplicate section, remove it
  and then recombine it (see Figure~\ref{fig explicit section
    delete}).  By reversing the section-deletion process we define
  \emph{duplication} of a section.

  We say that a set of animals, $\mathcal{A}$, is \emph{dense} if the
  set is closed under section deletion and duplication. \emph{ie} no
  animal outside the set can be produced by section deletion and / or
  duplication from a animal inside the set.
\end{defn}

The process of section-deletion and duplication leads to a partial
order on the set of animals.
\begin{defn}
  For any two animals $P, Q$, we write $P \nds Q$ if
  $P=Q$ or $P$ can be obtained from $Q$ by a sequence of
  section-deletions.  A \emph{section-minimal} animal, $P$, is a
  animal such that for all animals $Q$ with $Q \nds P$ we have
  $P=Q$.
\end{defn}
The above definition leads quite directly to the following lemma:
\begin{lemma}
  The binary relation $\nds$ is a partial order on the set of
  animals. Further every animal reduces to a unique section-minimal
  animal, and there are only a finite number of minimal animals with
  $n$ vertical bonds.
\end{lemma}

By considering the generating function of all animals that are equivalent (by some
sequence of section-deletions) to a given section-minimal animal, we
find that $H_n(x)$ may be written as the sum of simple rational
functions. Theorem~\ref{thm Hn nature} follows directly from
this. Further examination of the denominators of these functions gives
Theorem~\ref{thm poles sections}. Details are given in \cite{ADR_haru}.


\bibliographystyle{plain}
\bibliography{db_bib.bib}

\end{document}